\documentclass[12pt]{article}
\usepackage[a4paper, total={7in, 10in}]{geometry}
\usepackage[parfill]{parskip}
\usepackage{physics, tensor, float, subcaption}
\usepackage{graphicx}
\graphicspath{ {Plots/} }
\usepackage{jhep-mod}
\usepackage{bm}
\usepackage{soul}
\usepackage{amssymb,amsmath,amsthm}
\usepackage{mathrsfs}
\usepackage[utf8]{inputenc}
\usepackage{enumerate}
\usepackage{bigints}
\usepackage{xcolor}
\usepackage{appendix}
\usepackage{graphicx}
\usepackage{float}
\usepackage{tikz}
\usepackage{setspace}
\usepackage{cancel}
\usepackage{doi}
\definecolor{purple}{rgb}{1,0,1}
\definecolor{lime}{HTML}{A6CE39} 

\definecolor{lime}{HTML}{A6CE39}
\newcommand{\orcidicon}{%
	\begin{tikzpicture}
	\draw[lime, fill=lime] (0,0) 
		circle [radius=0.16] 
		node[white] {{\fontfamily{qag}\selectfont \tiny ID}};
	\draw[white, fill=white] (-0.0625,0.095) 
		circle [radius=0.007];
	\end{tikzpicture}
	\hspace{-5mm}
}
\newcommand\orcidMatt{{\href{https://orcid.org/0000-0003-1088-6485}{\orcidicon}}}

\renewcommand{\O}{\mathcal{O}}
\begin{document}

\title{\null\vspace{-50pt}
Primes in short intervals\\ assuming the Riemann hypothesis}

\author{
\Large
Matt Visser\!\orcidMatt\!
}
\affiliation{School of Mathematics and Statistics, Victoria University of Wellington, \\
\null\qquad PO Box 600, Wellington 6140, New Zealand.}
\emailAdd{matt.visser@sms.vuw.ac.nz}
\def\theta{\vartheta}
\def\O{{\mathcal{O}}}
\def\Li{{\mathrm{Li}}}

\abstract{
\vspace{1em}

In 2001 Baker, Harman, and Pintz proved,  for $x$ sufficiently large, 
the existence of at least one prime in the sub-linear interval 
$[x,x+x^{\frac{1}{2} +\frac{1}{40}}]$. But just how large is sufficiently large? 
In the current note I shall make the Baker--Harman--Pintz theorem effective, but at an extremely  high price. 
Contingent upon assuming validity of the Riemann hypothesis (RH),
I shall show that the Baker--Harman--Pintz theorem holds for $x \geq 3.5\times 10^{90}$.
Furthermore, again assuming RH, I shall also show that for any $0< m\leq15.831$ and $x \geq 1$ 
there is at least one prime in the interval  $[x,x+x^{\frac{1}{2} +\frac{1}{m}}]$, 
and that for  $M\geq 2.29$ and $y \geq 1$ 
there is at least one prime in the interval  $[y^M, (y+1)^M]$.

\bigskip

\bigskip
\noindent
{\sc Date:} 20 October 2026; \LaTeX-ed \today

}

\maketitle
\def\tr{{\mathrm{tr}}}
\def\diag{{\mathrm{diag}}}
\def\cof{{\mathrm{cof}}}
\def\pdet{{\mathrm{pdet}}}
\def\d{{\mathrm{d}}}
\parindent0pt
\parskip7pt
\def\Kerr{{\scriptscriptstyle{\mathrm{Kerr}}}}
\def\eos{{\scriptscriptstyle{\mathrm{eos}}}}
\section{Introduction}

Some 25 years ago Baker, Harman, and Pintz~\cite{BHP} proved existence of at least one prime in the sub-linear interval
\begin{equation}
\left[x, x+ x^{{1\over2}+{1\over 40}}\right] \qquad \hbox{for $x$ sufficiently large}.
\end{equation}
Baker, Harman, and Pintz commented that ``with enough work this result could be made explicit'', but as yet there has been no progress along these lines.
I shall make some progress on this point, but only at the very high price of assuming the Riemann hypothesis (RH). 
Assuming the RH, I will use the relatively recent 2019 Carneiro--Milinovich--Soundararajan bound on the prime gaps~\cite{gaps}, and combine this with the recent 2025 Sorenson--Webster range of validity for the upper part of the Oppermann conjecture~\cite{Oppermann}. 
I shall derive three theorems:\\
{\bf Theorem 1:}\\
Assuming the Riemann hypothesis, for $x\geq 3.5\times 10^{90}$ there is at least one prime in the interval 
$[x, x+ x^{{1\over2}+{1\over 40}}]$.\\
{\bf Theorem 2:}\\
Assuming the Riemann hypothesis, for $x\geq 1$  and $0< m\leq 15.831$ there is at least one prime in the interval 
$[x, x+ x^{{1\over2}+{1\over m}}]$.\\
{\bf Theorem 3:}\\
Assuming the Riemann hypothesis, for $y\geq 1 $  and $M\geq 2.29$ there is at least one prime in the interval 
$[y^M, (y + 1)^M]$.\\

\section{The Carneiro--Milinovich--Soundararajan bound }

One key result we will want to use is the 2019 Carneiro--Milinovich--Soundararajan bound on the prime gaps~\cite{gaps}.
Assuming RH one has
\begin{equation}
\label{E:gaps}
p_{n+1} < p_n + {22\over 25} \sqrt{p_n} \ln(p_n); \qquad (n \geq 2; \; p_n\geq 3).
\end{equation}
Now let $\pi(x)=n$, then $x \geq p_{\pi(x)} = p_n$ and $x< p_{\pi(x)+1} = p_{n+1}$. Thence for $n\geq 2$ we have
\begin{equation}
 x+{22\over 25} \sqrt{x} \ln(x) \geq  p_n {22\over 25} \sqrt{p_n} \ln(p_n) > p_{n+1} > x.
\end{equation}
Thence for $x \geq 3$ (assuming RH) there is at least one prime in the CMS interval
\begin{equation}
\left[ x,  x+{22\over 25} \sqrt{x} \ln(x) \right].
\end{equation}

\section{Comparing sub-linear intervals}
We now compare the generic BHP interval with the CMS interval
\begin{equation}
\left[ x,  x+x^{{1\over2}+{1\over m}} \right];
\qquad
\left[ x,  x+{22\over 25} \sqrt{x} \ln(x) \right].
\end{equation}
(For now assume $m>2$, to guarantee that the BHP interval is sub-linear.)
The two intervals are equal in size when
\begin{equation}
x^{1\over m} = {22\over 25}\ln(x); \qquad 
x = \exp\left(- m \; W_{-1}\left(-{25\over 22 \; m} \right) \right).
\end{equation}
Here $W_{-1}(z)$ is the $-1$ branch of the Lambert $W$ function.

Alternatively, solving for $m$ as a function of $x$, the two intervals are equal in size when
\begin{equation}
x^{1\over m} = {22\over 25}\ln(x); \qquad 
m = {\ln(x)\over \ln(\ln(x))+ \ln(22/25)}.
\end{equation}
The proofs of Theorems 1 and 2 are now straightforward.

\subsection{Proof of Theorem 1}

For the classic BHP theorem we have $m=40$ so crossover between the two intervals occurs at
\begin{equation}
x_{*} = \exp\left(- 40 \; W_{-1}\left(-{25\over 880} \right) \right)
=\exp\left(- 40 \; W_{-1}\left(-{5\over 176} \right) \right) = 3.482755638 \times10^{90}.
\end{equation}
Consequently for $x>x_*$ and $m= 40$ the BHP interval exceeds the CMS interval (which we already know contains primes).
So for $x> x_*$ the BHP interval contains primes. Rounding $3.482755638 \times10^{90}$ up to $3.5 \times10^{90}$ then yields Theorem 1. \hfill\hfill $\Box$

\subsection{Proof of Theorem 2}
We now compare the generic BHP interval with the interval appropriate for checking the upper part of Oppermann's conjecture: $\pi(n^2+n) > \pi(n^2)$, which is equivalent to the assertion $\pi(x+x^{1\over2}) > \pi(x)$. We note that for $m>0$:
\begin{equation}
\left[ x,  x+x^{{1\over2}+{1\over m}} \right]
\supset
\left[ x,  x+x^{1\over2} \right].
\end{equation}
Sorenson and Webster have recently verified~\cite{Oppermann} that the upper part of the Oppermann conjecture is valid for all 
$1\leq n \leq n_{**} = 7\times 10^{13}$, which corresponds to the interval $1\leq x_{**}=n_{**}^2 \leq 4.9\times 10^{27}$. 
So for arbitrary $m>2$ and $1\leq x \leq x_{**} = 4.9\times 10^{27}$ the generic BHP interval is guaranteed to contain at least one prime.

The case $m=2$ corresponds to the interval $[x,2x]$ which is covered by appealing to the Bertrand--Chebyshev theorem. 

The case $0<m<2$ corresponds to a super-linear interval larger than the  Bertrand--Chebyshev  interval $[x,2x]$.

Looking at the region above $x_{**}$, that is $x>x_{**}$, we want the generic BHP interval to exceed the CMS interval, which means we need to enforce
\begin{equation}
m <  {\ln(x_{**})\over \ln(\ln(x_{**}))+ \ln(22/25)} = 15.83179525.
\end{equation}
Rounding downwards, for $0<m\leq 15.831$ we are certain that the generic BHP interval will contain primes for all $x\geq 1$. 
This then yields Theorem 2. \hfill\hfill $\Box$

Let us now massage Theorem 2 to extract additional information.

\subsection{Proof of Theorem 3}
Let us set $x=y^M$ and choose $M=2m/(m-2)$ implying $m= 2M/(M-2)$.

Then it is easy to check that
\begin{equation}
M-1 = {2m\over m-2}-1 = {m+2\over m-2} = M {m+2\over2m} = M \left(\frac12 +\frac1m\right).
\end{equation}
Consequently under these conditions
\begin{equation}
[x, x+ x^{{1\over2}+{1\over m}}] = [y^M, y^M + y^{M-1}].
\end{equation}
Noting that $m= 15.83179525$ corresponds to $M= 2.289188780$, this now lets us re-write Theorem 2 as follows:\\
{\bf Corollary:}
Assuming the Riemann hypothesis, for $y\geq 1$  and $M\geq 2.289188780$ there is at least one prime in the interval 
$[y^M, y + y^{M-1}]$.
Rounding upwards, for $M\geq 2.29$  there is at least one prime in the interval 
$[y^M, y + y^{M-1}]$.\\
But since $M>1$ and $y\geq 1$ we certainly have
\begin{equation}
[y^M, y + y^{M-1}] \subset [y^M, y + M y^{M-1}] \subset [y^M, (y + 1)^M]
\end{equation}
Thence, assuming the Riemann hypothesis, for $y\geq 1$  and $M\geq2.29$ there is at least one prime in the interval 
$[y^M, (y + 1)^M]$.
This now yields Theorem 3. \hfill\hfill $\Box$

\section{Conclusions}

\enlargethispage{30pt}
Assuming the Riemann hypothesis is an extremely high price to pay --- but it does provide useful information.
\begin{description}
\item[Theorem 1:] This analysis  gives an explicit interval ($x \geq 3.5\times 10^{90}$) over which the original $m=40$ version of the Baker--Harman--Pintz result is guaranteed to hold. Given the strength of the input assumptions (RH) it is difficult to see how this might numerically be improved. Nevertheless, anyone wishing to derive unconditional effective bounds on the  Baker--Harman--Pintz result now has an idea of what region might plausibly be of interest.
\item[Theorem 2:] Improved validity regions for the Oppermann conjecture would allow one to increase the number $m=15.831$ appearing in Theorem 2. Doing so should be straightforward if tedious. Establishing complete validity of the Oppermann conjecture for all $x\geq 1$ would automatically allow $m$ to become infinite.
\item[Theorem 3:] Improved validity regions for the Oppermann conjecture would allow one to decrease the number $M=2.29$ appearing in Theorem 3. Doing so should be straightforward if tedious. Establishing complete validity of the Oppermann conjecture for all $x\geq 1$ would automatically  allow $M$ to take the value $2$.
\end{description}
Note that we have already implicitly shown that the exponent $M=3$ works for all $x\geq 1$ in Theorem 3 --- so assuming RH there are always primes between cubes. This is compatible with Lemma 5 of Caldwell--Cheng~\cite{Caldwell-Cheng}. See also the recent discussion in reference~\cite{almost-primes}. The value of the current Theorem 3 lies in the fact that we can push the result for the exponent $M=3$ down to $M=2.29$.

With weaker hypotheses, such as the 7/8 version of the quasi-Riemann hypothesis, considerably less can be said~\cite{QRH:7/8}. There is then no analogue of Theorem 1, while an analogue Theorem 2 holds for $m\geq 8/3$, and likewise an analogue of Theorem 3 holds for $M\geq 8$. 

\bigskip
\hrule\hrule\hrule

\clearpage

\bigskip
\hrule\hrule\hrule

\begin{thebibliography}{99}
\newcommand{\arXiv}[1]{arXiv:~{\href{https://arxiv.org/abs/#1}{\color{blue}#1}}}

\bibitem{BHP}
R. C. Baker, G. Harman, and J. Pintz,\\
``The Difference Between Consecutive Primes, II",\\
Proc. London Math. Soc. {\bf 83 \# 3} (2001) 532--562.

\bibitem{gaps}
Emanuel Carneiro, Micah B. Milinovich, Kannan Soundararajan,\\
``Fourier optimization and prime gaps'',\\
Comment. Math. Helv. {\bf 94 \# 3} (2019) 533--568,
[\arXiv{1708.04122} [math.NT]] \doi{10.4171/CMH/467}.

\bibitem{Oppermann}
Jonathan Sorenson, Jonathan Webster,\\
``An algorithm and computation to verify Legendre's Conjecture up to $7\times 10^{13}$'',\\
 Research in Number Theory {\bf 11} (2025) 4,
\url{https://doi.org/10.1007/s40993-024-00589-4}\\
(Preliminary version, only checking up to $3.33\times 10^{13}$, is available at \arXiv{2401.13753} [math.NT].)

\bibitem{Caldwell-Cheng}
C. K. Caldwell and Y. Cheng. 
``Determining Mills' Constant and a Note on Honaker's Problem''. 
J. Integer Seq. 8.05.4.1 (2005).
[\arXiv{1010.4883} [math.NT]].

\bibitem{almost-primes}
Daniel R. Johnston, Jonathan P. Sorenson, Simon N. Thomas, Jonathan E. Webster,
``Primes and almost primes between cubes'',
\arXiv{2601.15564} [math.NT].

\bibitem{QRH:7/8}
Matt Visser,
``Primes between 8th powers?'',
\arXiv{2610.20058}  math.NT

\end{thebibliography}
\end{document}